\documentclass{amsart}

\usepackage{amsmath}
\usepackage{amsthm}
\usepackage{amsfonts}
\usepackage{amssymb}
\usepackage{hyperref}

\theoremstyle{plain}
    \newtheorem{thm}{Theorem}[section]
    \newtheorem{prop}[thm]{Proposition}
    \newtheorem{lemma}[thm]{Lemma}
    \newtheorem{cor}[thm]{Corollary}

\theoremstyle{definition}

\theoremstyle{remark}

\numberwithin{equation}{section}

\newcommand{\rar}{\ensuremath{\rightarrow}}
\newcommand{\lrar}{\ensuremath{\longrightarrow}}

\newcommand{\Hom}{\textup{Hom}}
\newcommand{\StMod}{\textup{StMod}}

\newcommand{\Mod}{\textup{Mod}}

\newcommand{\uHom}{\underline{\Hom}}
\newcommand{\T}{\mathcal{T}}

\newcommand{\lstk}[1]{\stackrel{#1}{\longrightarrow}}
\newcommand{\Ext}{\textup{Ext}}
\newcommand{\Tate}{\widehat{H}}
\newcommand{\wOmega}{\widetilde{\Omega}}

\begin{document}

\title{Groups which do not admit ghosts}

\author{Sunil K. Chebolu}
\address{Department of Mathematics \\
University of Western Ontario \\
London, ON, Canada}
\email{schebolu@uwo.ca}

\author{J. Daniel Christensen}
\address{Department of Mathematics \\
University of Western Ontario \\
London, ON, Canada}
\email{jdc@uwo.ca}

\author{J\'{a}n Min\'{a}\v{c}}
\address{Department of Mathematics\\
University of Western Ontario\\
London, ON, Canada}
\email{minac@uwo.ca}

\keywords{Ghost map, stable module category,
derived category, Jennings' theorem, generating hypothesis.}
\subjclass[2000]{Primary 20C20, 20J06; Secondary 55P42}

\date{\today}

\commby{Paul Goerss}

\begin{abstract}
A ghost in the stable module category of a group $G$ is a map between
representations of $G$ that is invisible to Tate cohomology. We show that the
only non-trivial finite $p$-groups whose stable module categories have no non-trivial
ghosts are the cyclic groups $C_2$ and   $C_3$. We compare this to the
situation in the derived category of a commutative ring.
We also determine for which groups $G$
the second power of the Jacobson radical of $kG$ is stably isomorphic to a
suspension of $k$.
\end{abstract}

\maketitle

\section{Introduction}

Let $G$ be a $p$-group and let $k$ be a field of characteristic $p$.
A natural home for the cohomology of modular representations of $G$ is the stable module category
$\StMod(kG)$ of $G$. It  is the category obtained from the category of left $kG$-modules by killing the projectives.
It has the structure of a tensor triangulated category with the trivial representation $k$ as the unit object
and $\Omega$ as the loop (desuspension) functor.  The space of morphisms from $M$ to $N$ in $\StMod(kG)$ is denoted
$\uHom_{kG}(M,N)$ and consists of the $kG$-module homomorphisms modulo those that factor through a projective module.
Note that a map of
$kG$-modules that factors through a projective clearly induces the zero map in
Tate cohomology. A natural question is whether the converse is true. Namely,
if $f \colon M \rar N$ is a map of $kG$-modules such that the induced map in
Tate cohomology
\[\uHom_{kG}(\Omega^* k, M) \lrar \uHom_{kG}(\Omega^* k, N)\]
is trivial, then does $f$  factor through a projective module? Equivalently, is every such $f$ trivial
 in the stable module category? In this paper we investigate the groups $G$ for which the above question always has an
affirmative answer. Our main theorem states:

\begin{thm}
Let $G$ be a non-trivial finite $p$-group and let $k$ be a field of
characteristic $p$. All maps which are trivial in Tate cohomology
factor through a projective $kG$-module
if and only if $G$ is either $C_2$ or $C_3$.
\end{thm}

A map between $kG$-modules is said to be a \emph{ghost}
if the induced map in Tate cohomology $\uHom_{kG}(\Omega^* k, -)$ is trivial.
Using this terminology, the main theorem says that the only finite $p$-groups
for which all ghosts are trivial in the stable module category are $C_2$, $C_3$ and the trivial group.
It is interesting and surprising to note that the answer is so simple.

We now explain the strategy of the proof.
We construct a (weakly) universal ghost out of a given
$kG$-module $M$, and from this, we deduce that all ghosts out of $M$ vanish if
and only if $M$ is a retract of a direct sum of suspensions of $k$. In
particular, when $M$ is a finite-dimensional indecomposable module in
$\StMod(kG)$, all ghosts out of $M$ vanish if and only if $M \cong \Omega^i k$
in $\StMod(kG)$ for some $i$.  Therefore, we determine the finite $p$-groups $G$ which admit an
indecomposable $kG$-module that is not stably isomorphic to  any $\Omega^i k$.
A formula of Jennings~\cite{jennings}, which computes the nilpotency index of
the Jacobson radical of $kG$, plays an important role in this investigation.

The proof (sketched above) of our main result yields some interesting
additional results. For instance, 
Proposition~\ref{prop:J^2} characterises finite
$p$-groups $G$ for which the second power of the Jacobson radical of $kG$ is
isomorphic in the stable module category to a suspension of $k$.
In addition, the formal material 
in Section~\ref{sec:ghosts} has implications in other settings, as we
illustrate in Section~\ref{sec:derived}. 

We now explain how we were led to these results.  An old conjecture of Peter
Freyd~\cite{freydGH} in homotopy theory called the generating hypothesis (GH)
claims that a map $\phi \colon X \rar Y$ between finite spectra
that induces the zero map on stable homotopy groups is null-homotopic.
It is one of the most important  unsolved problems in stable homotopy theory.
In order to gain some insight into this deep problem, it is  natural to examine its analogues
in algebraic settings such as the derived category of a commutative ring
and the stable module category of  a finite group.
 The GH in the derived
 category $D(R)$ of a  ring $R$ is the statement that a map $\phi\colon X \rar
  Y$ between perfect complexes that induces the zero map in
 homology is chain-homotopic to the zero map.
Keir Lockridge~\cite{keir} showed that the GH holds in the derived category of a
commutative ring $R$ if and only if $R$ is a von Neumann regular ring.

In the stable module category, the GH is the statement that a map $\phi\colon M \rar N$ between
\emph{finite-dimensional} $kG$-modules is trivial in $\StMod(kG)$ if the induced map
in Tate cohomology is trivial.
This paper takes a first step towards the GH by studying the variant
in which the modules are not assumed to be finite-dimensional.
The theorem above implies that the GH is true for $C_2$ and $C_3$, but
does not answer the question for other groups, since it does not
guarantee the existence of ghosts between finite-dimensional modules.
Somewhat surprisingly, this variant of the GH turns out to be
equivalent to the GH.  Indeed,
motivated by the main result of this paper, we have shown in joint
work with Dave Benson~\cite{CCM3}
that the GH holds for a non-trivial finite $p$-group $G$ if and only if $G$ is either
$C_2$ or $C_3$.  This result can be used to deduce the main theorem stated above. However,
the methods used in the two papers are completely different. In this paper, we
use techniques inspired  from homotopy theory  and classical group theory,
while the techniques in~\cite{CCM3}
are more representation theoretic, relying heavily on the induction
and restriction. Moreover, we should emphasise
that the methods in~\cite{CCM3} do not give the additional results mentioned above.

We end the introduction by posing a riddle to the reader:
Which finite $p$-groups are like a finite product of fields?
We solve this riddle in the last section using a
result of Lockridge, which is  an analogue of our main theorem for the derived category of a commutative ring.
Thus the riddle sets a context for our main theorem, both in representation theory and  commutative algebra.

\section{Ghosts in triangulated
categories}\label{sec:ghosts}

Let $\T$ be a triangulated category which admits arbitrary coproducts and let
$S$ be a distinguished object. (If $\T$ is tensor triangulated,
we always take $S$ to be the unit object of $\T$.)  If $X$ and $Y$ are objects in $\T$,
then $[X, Y]_*$ will denote the graded abelian group of maps from $X$ to $Y$,
and $\pi_*(X)$ will stand for $[S, X]_*$. A map $\phi\colon M \rar N$ in $\T$
is a \emph{ghost} if the induced map
\[ \pi_*(M) \lstk{\phi_*} \pi_*(N) \]
is trivial. We say that $\phi \colon M \rar N$ is a \emph{universal ghost} if
$\phi$ is a ghost and if every ghost out of $M$ factors through
$\phi$. (Such a map should technically be called \emph{weakly universal} because we
don't assume the factorisation to be unique.) We begin by showing the
existence of a universal ghost out of any given object in $\T$.
This technique is well-known in homotopy theory, but we include
the details in this section to keep the paper self-contained.

Let $M$ be an object in $\T$. We assemble all the homogeneous elements of $\pi_*(M)$ into a map
\[\coprod_{\eta \, \in \, \pi_*(X)} \Sigma^{|\eta|}S \lrar M \]
from a coproduct of suspensions of the unit object,
where $|\eta|$ is the degree of $\eta$.
Completing this map to an exact triangle in $\T$, we get
\begin{equation} \label{eq:univ-ghost}
\coprod_{\eta \, \in \, \pi_*(X)} \Sigma^{|\eta|}S \lrar M \lstk{\Phi_M} U_M.
\end{equation}

\begin{prop}\label{prop:universal}
The map $\Phi_M\colon M \rar U_M$ is a universal ghost out of $M$.
\end{prop}

\begin{proof}
It is clear from its construction that $\Phi_M$ is a ghost. Now suppose $f\colon M \rar N$ is any ghost.
Then for each $\eta$ in $\pi_*(M)$, the composite
\[\Sigma^{|\eta|} S \lstk{\eta} M \lstk{f} N\]
is null. Therefore, so is the composite
\[\coprod_{\eta\,\in \, \pi_*(M)} \Sigma^{|\eta|} S \lrar M \lstk{f} N. \]
Now \eqref{eq:univ-ghost} is an exact triangle and therefore $f$ factors through $\Phi_M$.
\end{proof}

Proposition~\ref{prop:universal} essentially says
that ghosts form part of a projective class.
See~\cite{jdc} for more details.

Recall that an object $C$ in $\T$ is said to be \emph{compact} if the natural map
\[ \bigoplus_{\alpha} \ [C, X_{\alpha}] \lrar [C, \underset{\alpha}{\coprod}
 \ X_{\alpha}] \]
is an isomorphism for all set-indexed collections of objects $X_{\alpha}$ in
$\T$.  An object $X$ in $\T$ is said to be \emph{indecomposable} if $X \ne 0$  and a
decomposition  $X \cong A \amalg B$ in $\T$ implies
that either $A$ or $B$ is the zero object.

We say that $\T$ has the \emph{Krull-Schmidt property}
if the following two conditions hold:
\begin{itemize}
\item Each compact object in $\T$ can be decomposed uniquely into indecomposable objects.
\item The distinguished object $S$ in $\T$ is compact and indecomposable.
\end{itemize}
Our next proposition characterises the objects in $\T$ out of which
all ghosts vanish.

\begin{prop}\label{prop:trivial}
Let $\T$ be a triangulated category which admits arbitrary coproducts and
let $S$ be a distinguished object.
Then the following are equivalent for an object $M$ in $\T$:
\begin{enumerate}
\item All ghosts out of $M$ are trivial.
\item The universal ghost $\Phi_M\colon M \rar U_M$ is trivial.
\item $M$ is a retract of a coproduct of suspensions of $S$.
\end{enumerate}
Moreover, if $M$ is compact, then (3) is equivalent to:
\begin{enumerate}
\item[$(3')$] $M$ is a retract of a \emph{finite} coproduct of
suspensions of $S$.
\end{enumerate}
If $M$ is compact and $\T$ has the Krull-Schmidt property, then
(3) is equivalent to:
\begin{enumerate}
\item[$(3'')$] $M$ is a finite coproduct of suspensions of $S$.
\end{enumerate}
\end{prop}

\begin{proof}
(1) $\Rightarrow$ (2) is trivial, for $\Phi_M\colon M \rar U_M$ is itself a ghost. Now to see that
(2) implies (3), suppose the universal ghost  $\Phi_M$ is trivial. That
means that \eqref{eq:univ-ghost} is a split triangle. In particular, $M$ is a retract
of $\amalg \, \Sigma^i S$.  So there exists a map $M \lstk{j} \amalg \,\Sigma^i
\, S$ such that the composite
\[ M \lstk{j} \underset{ \pi_*(M)} {\coprod} \Sigma^i S  \lrar M \]
is the identity in $\T$.  Now if $M$ is compact, then  $j$ factors through a
finite coproduct. Therefore $M$ is a retract of a finite  coproduct of
suspensions of $S$. If $\T$ has the Krull-Schmidt property, then it  follows
that $M$ is a finite coproduct of suspensions of $S$. Finally (3)
$\Rightarrow$ (1) is clear.
\end{proof}

For reference we record the following corollaries which are immediate from  Proposition~\ref{prop:trivial}.

\begin{cor} \label{cor:ghostoutofindecomposable}
Let $\T$ be a triangulated category which admits arbitrary coproducts and
which has the Krull-Schmidt property, and let $S$ be a distinguished object. 
If $M$ is a compact indecomposable
object in $\T$ such that $M \ncong \Sigma^i S$ for any $i$, then there exists a
non-trivial ghost out of $M$.
\end{cor}

\begin{cor}  \label{cor:criterionfornoghosts}
Let $\T$ be a triangulated category which admits arbitrary coproducts
and let $S$ be a distinguished object. Every
ghost in  $\T$ is trivial if and only if $\T$ is the collection of retracts of
coproducts of suspensions of $S$.
\end{cor}

In the next section, we use these corollaries to determine when the stable
module category of a finite $p$-group has no non-trivial ghosts. In the
following section, we do the same for the derived category of a commutative
ring.

\section{Stable module categories}

We begin with some preliminaries. Let $G$ be a finite group and let $k$ be a
field. The stable module category $\StMod(kG)$ of $G$ is
the category obtained from the category of left $kG$-modules by killing off  the
projectives.  The space of morphisms from $M$ to
$N$ in $\StMod(kG)$ is denoted $\uHom_{kG}(M,N)$ and consists of the $kG$-module
homomorphisms modulo those that factor through a projective module.
Thus a map in the stable module category is trivial if and only
if it factors through a projective module. A key fact \cite{bencar-1992} is
that the Tate cohomology groups can be described as groups of morphisms in
$\StMod(kG)$: $\Tate^i(G,M) \cong \uHom(\Omega^i k, M)$. $\StMod(kG)$ has the
structure of a tensor triangulated category, where the trivial representation
$k$ is the unit object and $\Omega$ is the loop (desuspension) functor.
($\Omega M$ is defined to be the kernel of a projective cover of $M$. This is
well-defined in the stable module category.) We denote by $\wOmega^i(M)$
the projective-free part of $\Omega^i M$, which is a well-defined $kG$-module.
For more facts about $kG$-modules and $\StMod(kG)$, we refer the reader to
Carlson's excellent lecture notes \cite{carlson-modulesandgroupalgebras}.

From now on we work exclusively with finite $p$-groups and assume that
the characteristic of $k$ is $p$. We begin by proving the
easy direction of our main theorem.

\begin{prop} \label{prop:C_2C_3}
If $G$ is either $C_2$ or $C_3$, then $\StMod(kG)$ has no non-trivial ghosts.
\end{prop}

\begin{proof}
The group rings $kC_2 \cong k[x]/(x^2)$ and $kC_3 \cong k[x]/(x^3)$
are Artinian principal ideal rings. It is a fact \cite[p.~170]{injectivemodules} that every module
over an Artinian principal ideal ring is a direct sum of cyclic modules. We will use this fact to
show that every projective-free $kG$-module is a direct sum of suspensions of $k$. The result will
then follow from Corollary \ref{cor:criterionfornoghosts}.

First consider the group $C_2$. By the above fact, we know that every module over the ring $k[x]/(x^2)$
is a direct sum of copies of $k$ and $k[x]/(x^2)$. In particular, every projective-free $k[x]/(x^2)$-module  is a direct
sum of copies of $k$.  Now consider the group $C_3$.  In this
case, the above fact implies that every projective-free
$k[x]/(x^3)$-module is a direct sum of copies of $k$ and $k[x]/(x^2)$. But note that
\[\wOmega (k) := \ker\big(k[x]/(x^3) \twoheadrightarrow k\big)\cong  k[x]/(x^2).\]
In both cases we have shown that every projective-free $kG$-module is a direct
sum of suspensions of $k$. So we are done.
\end{proof}

We now collect some facts about finite $p$-groups that we need in the sequel.

\begin{lemma} \label{lem:indecomposabilitycriterion}
Let $G$ be a finite $p$-group and let $M$ be a finite-dimensional non-zero
$kG$-module. Then the invariant submodule $M^G$ of $M$ is non-zero. Thus there
is only one simple $kG$-module, namely the trivial module $k$. Moreover, if
$M^G$ is one-dimensional, then $M$ is indecomposable.
\end{lemma}

\begin{proof}
The proof of the first statement is an easy exercise; see \cite[3.14.1]{ben-1}. For the last statement, suppose
to the contrary that $M \cong A \oplus B $, with $A$ and $B$ non-zero. Then we have that $M^{P} \cong A^{P} \oplus  B^{P}$. This shows that
$M^{P}$ is at least two-dimensional, for by the first part of the lemma, both $A^{P}$ and $B^{P}$ are at least one-dimensional.
This contradiction completes the proof.
\end{proof}

\begin{lemma} \label{lem:ideals}
Let $G$ be a finite $p$-group and let $I$ be a non-trivial, proper ideal of
$kG$. Then $I$ is an indecomposable projective-free $kG$-module. In
particular, the powers $J^i(kG)$ which are non-zero are indecomposable
projective-free $kG$-modules.
\end{lemma}

\begin{proof}
We first show that $I$ is projective-free. If $I$ has a projective submodule,
then since projective modules over finite $p$-groups are free, that would mean that
$I$ should have $k$-dimension at least $|G|$, which is not possible since $I$ is
proper. To prove that $I$ is indecomposable, it suffices to show (by Lemma
\ref{lem:indecomposabilitycriterion}) that $I^G$ is one-dimensional. Note that
$I^G = I \bigcap \, (kG)^G$. It is easy to see that $(kG)^G$ is the
one-dimensional subspace generated by the norm element $\sum_{g \in G} g$. We
also know from Lemma~\ref{lem:indecomposabilitycriterion} that $I^G$ is
non-zero. Therefore $I^G$ is a one-dimensional submodule.
\end{proof}

\begin{lemma} \label{lem:dimofomegas}
Let $G$ be a finite $p$-group. Then, for all integers $i$, we have
\[ \dim(\wOmega^i k) \equiv (-1)^i \textup{ mod } |G|. \]
\end{lemma}

\begin{proof}
Recall that $\wOmega^1\, k$ is the kernel of the augmentation map, so we have a short exact sequence
\[ 0 \lrar \wOmega^1\, k \lrar kG \lrar k \lrar 0, \]
which tells us that $\dim (\wOmega^1\, k) \equiv  -1$ modulo $|G|$.
Inductively, it is clear from the short exact sequences
\[ 0 \lrar \wOmega^{i+1}\,k \lrar (kG)^t \lrar \wOmega^i\, k \lrar 0\]
that $\dim (\wOmega^i\,k) \equiv  (-1)^i$ modulo $|G|$ for $i \geq 0$.
(Here $(kG)^t$, for some $t$, is a minimal projective cover of $\wOmega^i \, k$.)
Also, since $\wOmega^i\,k \cong (\wOmega^{-i}\,k)^*$ in $\Mod(kG)$,
it follows that $\dim (\wOmega^i\,k) \equiv  (-1)^i$ modulo $|G|$ for
each integer $i$.
\end{proof}

We now introduce a formula of Jennings. Let $G$ be a finite $p$-group and let
$J(kG)$ be the Jacobson radical of $kG$. Since $kG$ is a local Artinian ring,
it follows that $J(kG)$ is nilpotent. So there exists a smallest integer $m$
such that $J(kG)^m = 0$. This integer will be called the \emph{nilpotency
index} of $J(kG)$, and it can be shown to be independent of the field $k$.
Very closely related  to the powers of the Jacobson radical are the dimension
subgroups of $G$, which we now define. The \emph{dimension subgroups} $F_i$ of
$G$ are defined by
\[ F_i := \{ g \in G : g -1 \in J^i(kG)\} \ \ \ \text{for} \;\; i \ge 1 .\]
These form a descending chain of normal subgroups in $G$
\[ F_1 \supseteq  F_2 \supseteq F_3 \supseteq \cdots \supseteq F_d \supseteq F_{d+1} ,\]
with $F_1 = G$ and $F_{d+1}$ trivial.
Define integers $e_i$ by $p^{e_i} = [F_i: F_{i+1}]$ for $1 \le i
\le d$. Then a formula due to Jennings~\cite{jennings} states that the
nilpotency index $m$ of $J(kG)$ is given by
\begin{equation} \label{for:nilpotency}
 m = 1 + (p - 1) \sum_{i = 1}^{ d}\, i\, e_i .
\end{equation}
Moreover, $e_1$ is the minimal number of generators for $G$. From the definition of the numbers $e_i$, it is clear that
$|G| = p^{\sum_i e_i}$.

\begin{prop} \label{prop:keyprop}
Let $G$ be a non-trivial finite $p$-group that is not isomorphic to $C_2$ or $C_3$. Then there exists a
finite-dimensional indecomposable projective-free $kG$-module  that is not isomorphic to
$\wOmega^i k$ for any $i$. In particular, there exists a non-trivial ghost in
$\StMod(kG)$.
\end{prop}

\begin{proof} It is well known that there are indecomposable projective-free modules over the Klein four group ($C_2 \oplus C_2$)
which are not of the form $\wOmega^i(k)$. In fact, every even-dimensional projective-free indecomposable $k(C_2 \oplus C_2)$-module
has this property.  Such modules are known to exist; see \cite[Thm.~4.3.3]{ben-1}, for instance.
Therefore, we can assume that $G$ is not any one of the groups
$C_2$, $C_3$ and $C_2 \oplus C_2$.

Consider the module $J^2(kG)$, which we denote by $J^2$ for brevity. By Lemma~\ref{lem:ideals}, we know that $J^2$ is an
indecomposable projective-free $kG$-module. Let $|G| = p^n$ for some positive integer $n$. Since the dimension of $\wOmega^i k$ is $+1$ or $-1$ modulo $|G| = p^n$
(see Lemma~\ref{lem:dimofomegas}), we will be done if we can show that the congruence class mod $p^n$ of $\dim(J^2)$  is different from $+1$ and $-1$. In fact, we
will show that when $G$ is not one of the above $3$ groups, then
\[ 1 < \dim(J^2) < p^n - 1 .\]
Note that $J^2 \subsetneq J$. For, otherwise, Nakayama's lemma would imply that $J = 0$, a contradiction. Therefore the second
inequality is clear. Now we establish the first inequality. We have two cases to consider here:
\medskip

\noindent
\textbf{Case 1:} Suppose $\dim J^2 = 0$. Then the nilpotency index $m$ of $J$ is $2$. So by Jennings' formula we have
\[ 2 = 1 + (p - 1) [e_1 + 2e_2 + \cdots + d e_d ].\]
This means $1 = (p - 1) [e_1 + 2e_2 + \cdots +  d e_d]$. Recall that $e_1$ is the minimal number of generators of $G$, so $e_1 > 0$.
Therefore, the last equation holds if and only if $p = 2$, $e_1 =1$ and $e_i = 0$ for $i \ge 2$. Since $|G| = p^{\sum_i e_i}$,
it follows that $J^2 = 0$ if and only if $G = C_2$. But $G \neq C_2$, by assumption. So this case cannot
arise.
\medskip

\noindent
\textbf{Case 2:} Suppose $\dim J^2 = 1$. That means $J^2 = k$, therefore
$J(J^2) = J(k) = 0$. So the nilpotency index $m$ of $J$ is $3$.
By Jennings' formula we have
\[ 3 = 1 + (p - 1) [e_1 + 2e_2 + \cdots + d e_d ].\]
This means $2 = (p - 1) [e_1 + 2e_2 + \cdots + d e_d ]$.  Here there are two
possibilities. Either $p = 3$, $e_1 = 1$, and $e_i = 0$ for all $i \ge 2$, or
$p=2$, $e_1 = 2$, and $e_i = 0$ for all $i \ge 2$. In the former, we have $|G| =
3$, so $G \cong C_3$, and in the latter, $|G| = 4$ and $G$ is
generated by $2$ elements, so $G \cong C_2 \oplus C_2$. But $G$ was assumed to
not be one of these groups, so this case also cannot arise.

Since both cases are ruled out, we have  proved the first inequality.

Finally, the last statement of the proposition follows from Corollary~\ref{cor:ghostoutofindecomposable}.
\end{proof}

The main step of the proof of Proposition~\ref{prop:keyprop} is essentially
the classification of $p$-groups with nilpotency index at most 3.
This is known to be an easy consequence of Jennings' formula,
but we have included a proof to keep the paper self-contained.

Combining Propositions~\ref{prop:C_2C_3} and~\ref{prop:keyprop},
we have proved our main theorem:

\begin{thm}
Let $G$ be a non-trivial finite $p$-group and let $k$ be a field of characteristic $p$. Then
$\StMod(kG)$ has no non-trivial ghosts if and only if $G$ is either
$C_2$ or $C_3$.
\end{thm}

We extract the following interesting result from the proof of Proposition~\ref{prop:keyprop}.

\begin{prop} \label{prop:J^2}
Let $G$ be a non-trivial finite $p$-group. Then $J^2(kG) \cong \wOmega^i k$ in
$\Mod(kG)$ for some $i$ if and only if $G$ is isomorphic to $C_3$ or $C_2
\oplus C_2$.
\end{prop}

\begin{proof}
The only if part follows from the proof of Proposition~\ref{prop:keyprop}. For
the converse, if $G$ is  $C_3$, then every finite-dimensional indecomposable
projective-free module (and hence $J^2(kG)$) is of the form $\wOmega^i k$, for
some $i$; see the proof of Proposition \ref{prop:C_2C_3}. If $G$ is $C_2
\oplus C_2$, then $J^2(kG) \cong k$. This completes the proof of the
proposition.
\end{proof}

\section{Derived categories}\label{sec:derived}

Let $R$ be a commutative ring and  let $D(R)$ be its (unbounded) derived
category. This has the structure of a tensor triangulated category with $R$
(viewed as a chain complex concentrated in degree zero) as the unit object.
Observe that a map $\phi\colon X \rar Y$ in $D(R)$ is a ghost if and only if
the induced map $H_*(\phi)\colon H_*(X) \rar H_*(Y)$ in homology is zero. The
natural question is to characterise commutative rings $R$ for which $D(R)$ has
no non-trivial ghosts. This has been done by Lockridge in~\cite{keir}. We
include a proof below for the reader's convenience and also to illustrate the
results in Section~\ref{sec:ghosts}.

\begin{thm} \cite{keir}
Let $R$ be a commutative ring. Then $D(R)$ has no non-trivial ghosts if and only if $R$
is a finite product of fields.
\end{thm}

\begin{proof}
Let $R$ be a finite product of fields, $F_1 \times F_2 \times \cdots \times
F_l$, say. Then every $R$-module splits naturally as a direct sum of modules
over the subrings $F_i$. It follows that $D(R)$ is  equivalent to $D(F_1)
\times D(F_2) \times \cdots \times D(F_l)$. Now note that the derived category
of a field $F$ is equivalent to the category of $\mathbb{Z}$-graded $F$-vector
spaces. From this it follows that every object in $D(R)$ is a retract of
direct sum of suspensions of $R$. Therefore, by Corollary
\ref{cor:criterionfornoghosts}, $D(R)$ does not have any non-trivial ghosts.
Conversely, suppose there are no non-trivial ghosts in $D(R)$. Then it is not
hard to see that for every pair of $R$-modules $M$ and $N$, we have
$\Ext_R^i(M, N) = 0$ for each $i > 0$. This implies that every $R$-module is
projective. Therefore $R$ is semi-simple; see \cite[Thm.~4.2.2]{wei}. Since
commutative semi-simple rings are precisely finite direct products of fields
(by the Artin-Wedderburn theorem), we are done.
\end{proof}

The answer to the riddle posed in the introduction should be now clear to the
reader. The only non-trivial finite $p$-groups that are like a finite product of fields
are $C_2$ and $C_3$.


\end{document}